\documentclass{mcom-l}

\newtheorem{theorem}{Theorem}[section]

\theoremstyle{definition}

\theoremstyle{remark}
\newtheorem{remark}{Remark}[section]

\numberwithin{equation}{section}

\usepackage{amsmath,amssymb,amsxtra,comment,graphicx,psfrag}
\usepackage{bm,mathrsfs}
\usepackage{mathtools}
\usepackage{color,xcolor}
\usepackage{enumerate,enumitem}
\usepackage{pdfsync}
\usepackage{hhline}
\usepackage{cite}

\usepackage{amsfonts}
\usepackage{algorithm}
\usepackage{algorithmicx}
\usepackage{algpseudocode}
\usepackage{yhmath}
\usepackage{qtree}
\usepackage{xcolor}
\usepackage{epsfig}
\usepackage{lineno,hyperref}

\usepackage[T1]{fontenc}
\usepackage[utf8]{inputenc}
\usepackage{booktabs}

\usepackage[toc]{appendix}
\numberwithin{equation}{section}
\usepackage{verbatim}\usepackage{graphicx}
\usepackage{geometry}
\usepackage[many]{tcolorbox}
\usepackage{caption}
\allowdisplaybreaks[1]{}

\newcommand{\be}{\begin{equation}}
\newcommand{\ee}{\end{equation}}
\newcommand{\ba}{\begin{array}}
\newcommand{\ea}{\end{array}}
\newcommand{\bea}{\begin{eqnarray*}}
\newcommand{\eea}{\end{eqnarray*}}
\newcommand{\bean}{\begin{eqnarray}}
\newcommand{\eean}{\end{eqnarray}}

\newcommand{\lc}{\mathrel{\raise2pt\hbox{${\mathop<\limits_{\raise1pt\hbox{\mbox{$\sim$}}}}$}}}
\newcommand{\gc}{\mathrel{\raise2pt\hbox{${\mathop>\limits_{\raise1pt\hbox{\mbox{$\sim$}}}}$}}}
\newcommand{\ec}{\mathrel{\raise1pt\hbox{${\mathop=\limits_{\raise2pt\hbox{\mbox{$\sim$}}}}$}}}

\newcommand{\nn}{\nonumber}

\begin{document}

\title{Energy diminishing implicit-explicit Runge--Kutta methods for gradient flows}

\author{Zhaohui Fu}
\address{Department of Mathematics, University of British Columbia, Vancouver, Canada}
\email{fuzhmath@gmail.com}

\author{Tao Tang}
\address{BNU-HKBU United International College, Zhuhai 519087, China, and Guangdong Provincial Key Laboratory of Interdisciplinary Research and Application for Data Science}
\curraddr{}
\email{tangt@sustech.edu.cn}

\author{Jiang Yang}
\address{Department of Mathematics, SUSTech International Center for Mathematics \& National Center for Applied Mathematics Shenzhen (NCAMS), Southern University of Science and Technology, Shenzhen 518055,China}
\curraddr{}
\email{yangj7@sustech.edu.cn}

\subjclass[2022]{}

\date{August 30, 2022}

\dedicatory{}


\begin{abstract}

This study focuses on the development and analysis of a group of high-order implicit-explicit (IMEX) Runge--Kutta (RK) methods that are suitable for discretizing gradient flows with nonlinearity that is Lipschitz continuous. We demonstrate that these IMEX-RK methods can preserve the original energy dissipation property without any restrictions on the time-step size, thanks to a stabilization technique. The stabilization constants are solely dependent on the minimal eigenvalues that result from the Butcher tables of the IMEX-RKs. Furthermore, we establish a simple framework that can determine whether an IMEX-RK scheme is capable of preserving the original energy dissipation property or not. We also present a heuristic convergence analysis based on the truncation errors. This is the first research to prove that a linear high-order single-step scheme can ensure the original energy stability unconditionally for general gradient flows. Additionally, we provide several high-order IMEX-RK schemes that satisfy the established framework. Notably, we discovered a new four-stage third-order IMEX-RK scheme that reduces energy. Finally, we provide numerical examples to demonstrate the stability and accuracy properties of the proposed methods.


\vskip .35cm
\end{abstract}

\maketitle

\section{Introduction}

Gradient flows are dynamics driven by a free energy. 
Many physical problems can be modeled  by PDEs that take the form of gradient flows, which are often derived from the second law of thermodynamics. To avoid a lengthy introduction of general gradient flow models, we will take phase field models as a typical  and concrete example. We also emphasize that all results obtained for phase field models  in this work are valid for general gradient flow models. Phase field models have emerged as a powerful approach for modeling and predicting mesoscale morphological and microstructural evolution in materials. They were originally derived for the microstructure evolution and phase transition, but have been recently extended to many other physical phenomena. Physically, it is well known that the energy functionals relevant to the phase field equations decay in time, which has been used as a useful criteria for designing numerical schemes. Consequently, great efforts have been devoted to the construction of efficient and accurate numerical methods preserving the energy decay property at the discrete level, most of them are of first- or second-order in time. The main purpose of this work is to study stability and convergence for a class of higher order methods for the time discretization of the phase field equations.
To demonstrate the main idea of our analysis, we consider a general class of phase-field models in the following form:

\begin{equation}\label{gf} 
    u_t = G(-D u+f(u)), \quad (x,t)\in \Omega\times [0,T],
\end{equation}
where $\Omega$ is a bounded domain in $\mathcal{R}^d$ ($d=1,2,3$), $T$ is a finite time, $G$ and $D$ are negative and dissipative operators, and the potential function $f$ satisfies Lipschitz conditions. For some well-known phase-field equations, we have 
\begin{itemize}
    \item for the Allen--Cahn equation, $G=-1, D=\epsilon^2\Delta, f(u)=u^3-u$;
        \begin{eqnarray}\label{AC}
        \frac{\partial u}{\partial t} =\epsilon^2 \Delta u-f(u);
        \end{eqnarray}
    \item for the Cahn--Hilliard equation, $G=\Delta, D=\epsilon^2\Delta, f(u)=u^3-u$;
         \begin{eqnarray}\label{CH}
         \frac{\partial u}{\partial t} =\Delta (-\epsilon^2 \Delta u+f(u));
         \end{eqnarray}
    \item for the MBE (molecular beam epitaxy) model without slope selection or also known as the thin film model, $G=-1, D=-\epsilon^2\Delta^2, f(u)=-\nabla \cdot(\frac{\nabla u}{1+|\nabla u|^2})$
         \begin{eqnarray}\label{MBE}
         \frac{\partial u}{\partial t}=-\epsilon^2 \Delta^2 u+\nabla \cdot(\frac{\nabla u}{1+|\nabla u|^2}).
         \end{eqnarray}
\end{itemize}

In the Allen--Cahn model (\ref{AC}) and the Cahn--Hilliard model (\ref{CH}), $u$ represents the concentration of one of two components in the alloy or certain phase field, and the parameter $\epsilon$ measures the interfacial width, which is small compared to the characteristic length of the laboratory scale. In the thin film model (\ref{MBE}), the function $u$ is the scaled height of epitaxial growth thin films in a co-moving frame. It is well-known that these models satisfy the energy dissipation law, since all of them can be viewed as the gradient flows with the following energy functionals respectively:
\begin{equation}\label{EAC}
    \mathcal{E}(u)=\int_\Omega \left(\frac{\epsilon^2}{2} | \nabla u|^2 +F(u)\right) dx
\end{equation}
in $L^2$ for the Allen--Cahn equation, in $H^{-1}$ for the Cahn--Hilliard equation, and
\begin{equation}\label{EMBE}
    \mathcal{E}(u)=\int_\Omega \left(\frac{\epsilon^2}{2} | \Delta u|^2 +F(\nabla u)\right) dx
\end{equation}
for the thin film model, where $F(\nabla u)=\frac{1}{2} \ln(|\nabla u|^2+1)$.

Due to the perturbed (i.e., the scaling coefficient $\epsilon^2\ll 1$) Laplacian or biharmonic operators and strong nonlinearities involved in the equations, it is difficult to design an efficient and accurate time discretization scheme which resolve dynamics and steady states of the corresponding phase-field models. Moreover, another challenging issue for numerical approximations is to preserve the energy dissipation law which intrinsically holds for all these models. Numerical evidence has shown that non-physical oscillations may happen when the energy stability is violated. Therefore, a satisfactory numerical strategy needs to balance accuracy, efficiency and nonlinear stability of the solution.

There have been a wide range of studies for the construction of various numerical schemes preserving the energy dissipation law at a discrete level. Some popular and significant implicit time stepping work includes convex splitting methods \cite{epitaxy3,conv2} and Crank-Nicolson type scheme \cite{CHCN1,CHCN2}. The deficiency of these methods is the expense of solving a nonlinear system of equations at each time step. While in contrast to fully implicit schemes, implicit-explicit (also named semi-implicit) methods treats the nonlinear term explicitly and the linear term implicitly, and to solve these linearly implicit problems, it only requires solving a linear system of equations at every time step. Such methods may date back to the work of Chen and Shen \cite{IMEXpf1} in the phase-field context, and so far there have been developed many techniques and skills to design such schemes, see, e.g., \cite{SAVRK,ACCH,splitting1,robust,ATA,opspltLi,yang2015SDC,xuyan,Wangxm,exRKAC}. Based on the idea of the invariant energy quadratization (IEQ) method (see, e.g., \cite{IEQ1,IEQ2}), Shen et al. \cite{robust,SAV2} proposed the scalar auxiliary variable (SAV) method, which can easily render the unconditional energy decay property. Some modified SAV methods are developed recently \cite{SAVxu,SAVcheng,SAVqiao,TWY-JSC22}. However, the energy considered in these methods is modified from the original energy. In another direction, exponential time differencing (ETD) methods for the Allen--Cahn equation and other semilinear parabolic equations have attracted much recent attention. Du et al. \cite{ETD1} has shown that ETD and ETDRK2 schemes unconditionally preserve maximum bound property (MBP) and energy stability (but not the dissipation law). In particular, \cite{fuyangetd} establishes the original energy stability for ETDRK2. For the thin film model (or MBE model), some interesting results concerning stability analysis and error estimates for the ETD schemes are obtained in \cite{ETDstiff,ETDanal,epitaxial1,epitaxy4,epxywithornot}. Moreover, it is shown in \cite{fimRk} that fully implicit Runge-Kutta methods can decrease the energy of gradient systems, but the existence and uniqueness of the solution are still open problems. Another class of implicit Runge--Kutta methods for phase-field models is based on convex splitting approach, which demonstrates some good stability properties \cite{convexsplitphasefield}.
 
In this work, we prove that the energy dissipation law is unconditionally satisfied by a class of high-order IMEX-RK schemes for gradient flows with Lipschitz nonlinearity. We make use of a stabilization technique and the relationship between the stabilization coefficients and the time step is established. As long as the conditions for the RK method are satisfied, the schemes under consideration satisfy the energy stability unconditionally. 
Thus, the IMEX-RK method becomes the first provable high-order linear one-step scheme which unconditionally decreases the {\bf original energy} for the gradient flows. Furthermore, an error estimate based on the truncation error is provided.

This paper is organized as follows. In Section \ref{sec2}, we introduce some preliminaries concerning splitting techniques and IMEX-RK schemes. We present and prove our main stability theorems related to the Allen--Cahn and the Cahn--Hilliard equation in Section \ref{sec3}. The error analysis is provided in Section \ref{sec4}. Section \ref{sec5} offers some IMEX-RK examples to illustrate the requirements of the theorems, followed by several numerical examples in Section \ref{sec6}. Some concluding remarks are given in the final section.

\section{Preliminaries}\label{sec2}

\subsection{A stabilization technique}
In this section we introduce some preliminaries concerning operator splitting and implicit-explicit (IMEX) Runge-Kutta (RK) schemes. We consider the general form of gradient flows
\begin{equation} 
    u_t = G(- D u+f(u)), \quad (x,t)\in \Omega\times [0,T]
\end{equation}
where $G$ and $D$ are negative-definite operators and the function $f$ is Lipschitz continuous with constant $L$, i.e. given any $u, v$,
\begin{equation}
\begin{aligned}
        F(u)-F(v) &\leq (u-v)f(v)+\frac{L}{2} (u-v)^2
\end{aligned}
\end{equation}
for the Allen--Cahn and Cahn--Hilliard equations and 
\begin{equation}
\begin{aligned}
        F(\nabla u)-F(\nabla v) &\leq \partial_{\nabla v}F(\nabla v)\cdot(\nabla u- \nabla v)+\frac{L}{2} |\nabla u-\nabla v|^2
\end{aligned}
\end{equation}
for the MBE model without slope selection. The Lipschitz requirement is usually satisfied by either assuming or proving the solution is $L^\infty$ bounded or the function $f(u)$ is Lipschitz continuous for all $u\in \mathcal{R}$. In this paper, $(\cdot,\cdot)$ represents the spatial inner product equipped to the $L^2$ space.
\begin{remark}
    For simplicity we just assume the nonlinear term satisfies the Lipschitz condition, while in fact it only needs to be locally Lipschitz continuous in a strip along the exact solution (see assumptions in \cite{ExpoIntg}). Notice that for the Cahn–Hilliard equation and some other phase field models, the potential function $f(u)$ does not satisfy the Lipschitz assumption. However, in practical computations, the solution $u$ is often bounded and if we truncate $F(u)$ to quadratic growth for sufficiently large $|u|$, then the numerical solution is exactly the same as the one without the truncation, see \cite{stabsizeforsemiFrCH,2ndsemi2DCH,convganalfornonlocalCH}. It can be easily verified that the truncated potential function satisfies the Lipschitz assumption.
\end{remark}
For simplicity, periodic boundary conditions or homogeneous Neumann boundary conditions are imposed. Besides, for the Cahn--Hilliard equations, the zero mean initial data $u|_{t=0}=u_0$ is considered to guarantee that $\Delta^{-1}$ and other related operators are well-defined. The energy of the gradient flow is defined by 
\begin{equation}
    \mathcal{E}(u)=\int_\Omega\left( \frac{1 }{2}  |D_{1/2} u |^2 +  F(u) \right)\ dx,
\end{equation}
where $D_{1/2}D_{1/2}^*=D_{1/2}^* D_{1/2}=-D $ and $F'=f$. Consider the splitting of the energy $\mathcal{E}(u)=\mathcal{E}_l-\mathcal{E}_n$ with
\begin{eqnarray}
    && \mathcal{E}_l(u)=\int_\Omega\left(  \frac{1 }{2}  |D_{1/2} u |^2 + \frac{\alpha}{2}  |D_{1/2} u |^2 + \frac{\beta}{2} |u|^2 \right)\ dx, \\
    && \mathcal{E}_n(u)=\int_\Omega\left(  -   F(u) + \frac{\alpha }{2}  |D_{1/2} u |^2 + \frac{\beta}{2} |u|^2  \right)\ dx,
\end{eqnarray}
In the computation, $\mathcal{E}_l$ is treated implicitly and $\mathcal{E}_n$ is treated explicitly, which leads to an implicit scheme but still linear to solve. From this perspective, the gradient flow can be equivalently written as
\begin{equation}\label{sgf} 
    u_t=G(D_s u-f_s(u)),
\end{equation}
where $D_s=-(1+\alpha)D+\beta I$ and $f_s=-f-\alpha D+\beta I$.


\subsection{Implicit-Explicit Runge-Kutta schemes}
For the linear term $G D_s u$ in the gradient flow (\ref{sgf}), we consider an $s$-stage diagonally implicit Runge-Kutta (DIRK) scheme with coefficients $A=\left(a_{ij}\right)_{s\times s}\in \mathcal{R}^{s\times s},c,b\in\mathcal{R}^{s}$, in the usual Butcher notation. For the nonlinear term $G f_s (u)$ in (\ref{sgf}), we make use of an $s$-stage explicit scheme with the same abscissae $\hat{c}=c$ and coefficient $\hat{A}=\left(\hat{a}_{ij}\right)_{s\times s}\in \mathcal{R}^{s\times s},\hat{b}\in\mathcal{R}^{s}$. Thus, the IMEX-RK scheme can be determined by the following Butcher notation

\begin{equation}\label{IMEXRK}
\begin{array}{c|ccccc}
    0   & 0 & 0      &...& ... & 0 \\
    c_1 & 0 & a_{11} & 0 & ... & 0 \\
    c_2 & 0 & a_{21} & a_{22} & ... & 0 \\
    ... & 0 & ... & ... & ... & ...\\
    c_s & 0 & a_{s1} & a_{s2} & ... & a_{ss}\\
    \hline
        & 0 & b_1 & b_2 & ... & b_s
\end{array} \quad 
\begin{array}{c|ccccc}
        0     &     0    & 0 &  0   & ... & 0 \\
    \hat{c}_1 & \hat{a}_{11} & 0 & 0 & ... & 0 \\
    \hat{c}_2 & \hat{a}_{21} & \hat{a}_{22} & 0 & ... & 0 \\
    ... & ... & ... & ... & ... & 0\\
    \hat{c}_s & \hat{a}_{s1} & \hat{a}_{s2} & ... & \hat{a}_{ss} & 0\\
    \hline
        & \hat{b}_1 & \hat{b}_2 & ... & \hat{b}_s & 0
\end{array}.
\end{equation}
where $c_i=\hat{c}_i=\sum_j a_{ij}=\sum_j \hat{a}_{ij}.$ Here we only consider such special IMEX-RK schemes: $b_j=a_{sj},\hat{b}_j=\hat{a}_{sj},j=1,...,s$, indicating the implicit scheme is stiffly accurate. Besides, we require that the coefficient matrix of the implicit scheme begins with a zero column, i.e. we are considering type ARS (see \cite{Ascher}). In addition, we also require that the matrix $\hat{A}$ is invertible.

Consider a model equation:
\begin{equation}
    u_t=\mathcal{L}u+N(u),
\end{equation}
where $\mathcal{L}$ represents a linear operator and $N$ indicates a nonlinear one. Applying the IMEX-RK (\ref{IMEXRK}) to the model equation, we derive the following system of equations: (to solve $u_{n+1}$ starting from $u_n$)
\begin{equation}\label{IMEXRKsys}
    \left\{\begin{aligned}
    &v_0=u_n,\\
    &v_i=v_0+\tau \left( \sum_{j=1}^i a_{ij} \mathcal{L} v_j +\sum_{j=1}^i \hat{a}_{ij} N(v_{j-1})  \right), \quad 1\leq i \leq s, \\
    &u_{n+1}=v_s.
    \end{aligned}\right.
\end{equation}
There have been some existing work related to such IMEX-RK schemes, such as asymtotic behavior, error analysis and studies of different kinds of stability, etc. See, e.g. \cite{Ascher,betaConv,errestIMEXRK,DIRKMOL}.

\section{Energy decreasing property}\label{sec3}

In this section we present our main theorem and its proof.
\subsection{Main Theorem}

\begin{theorem} \label{maintheorem}
The IMEX-RK scheme (\ref{IMEXRKsys}) unconditionally decreases the energy of the phase field model (\ref{sgf}) if the following three matrices are positive-definite:
\begin{eqnarray}
    && H_0= (\hat{A})^{-1} E_{L},\label{H0}\\
    && H_1(\beta)=\beta Q - \frac{L}{2} I, \label{H1}\\
    && H_2(\alpha)=\alpha Q+(\hat{A})^{-1} A E_{L}-\frac{1}{2} E.\label{H2}
\end{eqnarray}
Here we say a matrix $M$ is positive-definite when $\frac{1}{2} (M+M^T)$ is positive-definite. In (\ref{H0} - \ref{H2}) $\alpha \text{ and } \beta$ are stabilizer constants, $E=\mathbf{1}_{s\times s}$, $E_{L}=\left(\mathbf{1}_{i\geq j}\right)_{s\times s}$ represents the lower triangular matrix of elements $1$, $I=\left(\mathbf{1}_{i=j}\right)_{s\times s}$ represents the identity matrix, and the determinant matrix $Q$ is defined by
\begin{equation}
    Q=\left((\hat{A})^{-1} A-I\right)E_{L}+I.
\end{equation}
In the case that both $Q$ and $H_0$ are positive-definite, if $\alpha$ and $\beta$ are sufficiently large in the sense that
\begin{equation}\label{alpbet}
\begin{aligned}
    &\beta  \geq \beta_0 := \frac{L}{2 \lambda_{\min}(Q)},\\
    &\alpha \geq \alpha_0 :=- \frac{\lambda_{\min}((\hat{A})^{-1} A E_{L} -\frac{1}{2}E)}{\lambda_{\min}(Q)},
\end{aligned}
\end{equation}
where $\lambda_{\min}$ represents the smallest eigenvalue, then $H_0, H_1$ and $H_2$ are all positive-definite. Consequently, the IMEX-RK scheme (\ref{IMEXRKsys}) satisfies the energy diminishing property unconditionally.
\end{theorem}

\begin{proof}

The system of the IMEX-RK schemes (\ref{IMEXRKsys}) can formally be written as
\begin{equation}\label{IMEXRKsysalter}
    \begin{pmatrix}
    v_1 \\ v_2 \\ \vdots \\ v_s
    \end{pmatrix}
    =
    \begin{pmatrix}
    v_0 \\ v_0 \\ \vdots \\ v_0
    \end{pmatrix}
    +\tau\left(A \begin{pmatrix}
    \mathcal{L}v_1 \\ \mathcal{L}v_2 \\ \vdots \\ \mathcal{L}v_s
    \end{pmatrix}+\hat{A} \begin{pmatrix}
    N(v_0) \\ N(v_1) \\ \vdots \\ N(v_{s-1})
    \end{pmatrix}    \right).
\end{equation}
Thus, we have
\begin{equation}
    \frac{1}{\tau} \begin{pmatrix}
    v_1-v_0 \\ v_2-v_0 \\ \vdots \\ v_s-v_0
    \end{pmatrix}
    =A \begin{pmatrix}
    \mathcal{L}v_1 \\ \mathcal{L}v_2 \\ \vdots \\ \mathcal{L}v_s
    \end{pmatrix}
    +\hat{A} \begin{pmatrix}
    N(v_0) \\ N(v_1) \\ \vdots \\ N(v_s)
    \end{pmatrix}.
\end{equation}
For the phase field model, we have $ \mathcal{L}v=G D_s v,\;N(v)=-G f_s(v)$, where $D_s=-(1+\alpha)D+\beta I$ and $f_s=-f-\alpha D+\beta I$. So we can derive the following equation
\begin{equation}\label{bigequation}
\begin{aligned}
    \frac{1}{\tau} \begin{pmatrix}
    v_1-v_0 \\ v_2-v_0 \\ \vdots \\ v_s-v_0
    \end{pmatrix}
    =&-A \begin{pmatrix}
    GDv_1 \\ GDv_2 \\ \vdots \\ GDv_s
    \end{pmatrix}
    +\hat{A} \begin{pmatrix}
    Gf(v_0) \\ Gf(v_1) \\ \vdots \\ Gf(v_{s-1})
    \end{pmatrix}  
     -\alpha \left[A \begin{pmatrix}
    GDv_1 \\ GDv_2 \\ \vdots \\GDv_s
    \end{pmatrix}  -\hat{A}\begin{pmatrix}
    GDv_0 \\ GDv_1 \\ \vdots \\GDv_{s-1}
    \end{pmatrix}\right] \\
    &+\beta \left[A \begin{pmatrix}
    Gv_1 \\ Gv_2 \\ \vdots \\ Gv_s
    \end{pmatrix} -\hat{A}\begin{pmatrix}
    Gv_0 \\ Gv_1 \\ \vdots \\ Gv_{s-1}
    \end{pmatrix}\right].
\end{aligned}
\end{equation}

Below we will reformulate the system. For simplicity, we denote the inverse matrix of $\hat{A}$ as $\hat{B}\in \mathcal{R}^{s\times s}$,
\begin{equation}
    \hat{A}\hat{B}=\hat{B}\hat{A}=I.
\end{equation}

In (\ref{rslt1} - \ref{rslt3}) below, we list some simple but useful results to help the simplification procedure.

\begin{equation}\label{rslt1}
\begin{pmatrix}
    v_1-v_0 \\ v_2-v_0 \\ \vdots \\ v_s-v_0
    \end{pmatrix}= E_L    \begin{pmatrix}
    v_1-v_0 \\ v_2-v_1 \\ \vdots \\ v_s-v_{s-1}
    \end{pmatrix}.
\end{equation}

\begin{equation}\label{rslt2}
     (A-\hat{A})\mathbf{1}_{s\times1} = c-\hat{c}= 0,
\end{equation}
where $\mathbf{1}_{s\times 1}=(1,1,...,1)^T \in \mathcal{R}^{s\times1}$. The equation (\ref{rslt1}) shows the matrices $A$ and $\hat{A}$ share the same eigenvector $\mathbf{1}_{s\times 1}$ which is useful in simplification. Taking the second term on the right hand side of (\ref{bigequation}) as an example, 
\begin{eqnarray}\label{rslt2}
    A \begin{pmatrix}
    GDv_1 \\ GDv_2 \\ \vdots \\ GDv_s
    \end{pmatrix}&=&A \begin{pmatrix}
    GDv_1 \\ GDv_2 \\ \vdots \\ GDv_s
    \end{pmatrix}-(A-\hat{A})\begin{pmatrix}
    GDv_0 \\ GDv_0 \\ \vdots \\ GDv_0
    \end{pmatrix} \\
    &=&A\begin{pmatrix}
    GDv_1-GDv_0 \\ GDv_2-GDv_0 \\ \vdots \\ GDv_s-GDv_0
    \end{pmatrix}+\hat{A}\begin{pmatrix}
    GDv_0 \\ GDv_0 \\ \vdots \\ GDv_0
    \end{pmatrix}\nonumber \\
    &=&A E_L \begin{pmatrix}
    GD(v_1-v_0) \\ GD(v_2-v_1) \\ \vdots \\ GD(v_s-v_{s-1})
    \end{pmatrix}+\hat{A}\begin{pmatrix}
    GDv_0 \\ GDv_0 \\ \vdots \\ GDv_0
    \end{pmatrix}.\nonumber
\end{eqnarray}

We also notice that
\begin{eqnarray}\label{rslt3}
    &&A \begin{pmatrix}
    Gv_1 \\ Gv_2 \\ \vdots \\ Gv_s
    \end{pmatrix} -\hat{A}\begin{pmatrix}
    Gv_0 \\ Gv_1 \\ \vdots \\ Gv_{s-1}
    \end{pmatrix} \\
    &=&\left[A \begin{pmatrix}
    Gv_1 \\ Gv_2 \\ \vdots \\ Gv_s
    \end{pmatrix} -\hat{A}\begin{pmatrix}
    Gv_0 \\ Gv_1 \\ \vdots \\ Gv_{s-1}
    \end{pmatrix}\right]
    -\hat{A}\begin{pmatrix}
    Gv_1 \\ Gv_2 \\ \vdots \\ Gv_s
    \end{pmatrix}
    +\hat{A}\begin{pmatrix}
    Gv_1 \\ Gv_2 \\ \vdots \\ Gv_s
    \end{pmatrix} \nonumber \\
    &=&(A-\hat{A})\begin{pmatrix}
    Gv_1 \\ Gv_2 \\ \vdots \\ Gv_s
    \end{pmatrix}
    +\hat{A} \begin{pmatrix}
    Gv_1-Gv_0 \\ Gv_2-Gv_1 \\ \vdots \\ Gv_s-Gv_{s-1}
    \end{pmatrix}\nonumber \\
    &=&(A-\hat{A})\left[\begin{pmatrix}
    Gv_1 \\ Gv_2 \\ \vdots \\ Gv_s
    \end{pmatrix}-\begin{pmatrix}
    Gv_0 \\ Gv_0 \\ \vdots \\ Gv_0
    \end{pmatrix}\right]
    +\hat{A} \begin{pmatrix}
    Gv_1-Gv_0 \\ Gv_2-Gv_1 \\ \vdots \\ Gv_s-Gv_{s-1}
    \end{pmatrix} \nonumber\\
    &=&(A-\hat{A}) E_L (Gw)+\hat{A} (Gw)=\hat{A} Q (Gw), \nonumber
\end{eqnarray}

where for simplicity we denote \[w=(v_1-v_0,...,v_s-v_{s-1})^T,\quad\quad\quad\quad\] \[Q=\left((\hat{A})^{-1} A-I\right)E_L+I,\quad\quad\quad\quad\quad\] and \[(Pw)=(P(v_1-v_0),...,P(v_s-v_{s-1}))^T\] for all operators $P$.

Therefore, using all above equations and recall that $\hat{B}=\hat{A}^{-1}$, (\ref{bigequation}) can be reformulated in the following way by applying $G^{-1}$ to both sides
\begin{equation}
    \frac{1}{\tau} E_L (G^{-1}w)=-A E_L (Dw)-\hat{A} \begin{pmatrix}
    Dv_0 \\ Dv_0 \\ \vdots \\ Dv_0
    \end{pmatrix}+\hat{A} \begin{pmatrix}
    f(v_0) \\ f(v_1) \\ \vdots \\ f(v_{s-1})
    \end{pmatrix} 
    -\alpha \hat{A} Q (Dw)+\beta \hat{A} Q (w).
\end{equation}
Consequently, one can isolate the nonlinear term as
\begin{equation}\label{nonlinearterm}
    \begin{pmatrix}
    f(v_0) \\ f(v_1) \\ \vdots \\ f(v_{s-1})
    \end{pmatrix} = {\frac{1}{\tau}} \hat{B}  E_L (G^{-1}w)
    +\hat{B} A E_L (Dw) + \begin{pmatrix}
    Dv_0 \\ Dv_0 \\ \vdots \\ Dv_0
    \end{pmatrix}
    +\alpha Q (Dw)-\beta Q (w).
\end{equation}

Next we focus on the difference of the energy, which reads
\begin{equation}\label{diffE}
    \mathcal{E}_{n+1}-\mathcal{E}_n= -\frac{1}{2} \left((u_{n+1},Du_{n+1}) -(u_n,Du_n)\right)+(F_{n+1}-F_n,\mathbf{1}),
\end{equation}
where $(u,v)$ represents the inner product in space and $\mathbf{1}$ is $\mathbf{1}_{s\times1}$ for simplicity. Notice that given any $u,v$, we have
\begin{equation}
    F(u)-F(v)\leq f(v) (u-v)+\frac{L}{2}  (u-v)^2,
\end{equation}
where $L$ is the Lipschitz constant of the function $f$.
Thus,
\begin{eqnarray}\label{Lipineqrslt}
        &&(F_{n+1}-F_n,\mathbf{1}) = \sum_{i=0}^{s-1} \ (F(v_{i+1})-F(v_i),\mathbf{1})\\
    &\leq& \begin{pmatrix}
    v_1-v_0, v_2-v_1,\cdots,v_s-v_{s-1}
    \end{pmatrix} \begin{pmatrix}
    f(v_0) \\ f(v_1) \\ \vdots \\ f(v_{s-1})
    \end{pmatrix}+\frac{L}{2}\begin{pmatrix}
    v_1-v_0, v_2-v_1,\cdots,v_s-v_{s-1} 
    \end{pmatrix}^2\nonumber \\
    &=& w^T  \begin{pmatrix}
    f(v_0) \\ f(v_1) \\ \vdots \\ f(v_{s-1})
    \end{pmatrix} + \frac{L}{2} w^T w,\nonumber
\end{eqnarray}
where the product for vector functions here is defined by $a^T b=\sum_{i=1}^n (a_i,b_i)$. Substituting $\begin{pmatrix}
    f(v_0), f(v_1), ...,f(v_{s-1})
    \end{pmatrix}^T$ to (\ref{Lipineqrslt}) by using (\ref{nonlinearterm}) gives
\begin{eqnarray}
        &&(F_{n+1}-F_n,\mathbf{1})\\
    &\leq &\frac{L}{2} w^T w + w^T \left(\frac{1}{\tau} \hat{B} E_L (G^{-1}w)
    +\hat{B} A  E_L (Dw) +\begin{pmatrix}
    Dv_0 \\ Dv_0 \\ \vdots \\ Dv_0
    \end{pmatrix}\right)
    +\alpha w^T  Q (Dw)-\beta w^T Q (w) \nonumber
\end{eqnarray}
Notice that
\begin{equation}\label{simple1rslt}
    w^T\begin{pmatrix}
    Dv_0 \\ Dv_0 \\ \vdots \\ Dv_0
    \end{pmatrix} =(v_s-v_0, Dv_0).
\end{equation}

Combining the above result (\ref{simple1rslt}) with the quadratic term in the energy difference (\ref{diffE}) gives
\begin{eqnarray}
    &&-\frac{1}{2} \left( (v_s,Dv_s)-(v_0,Dv_0) \right)+(v_s-v_0, Dv_0)\\
&=&-\frac{1}{2}(v_s,Dv_s)    -\frac{1}{2}(v_0,Dv_0)+(v_s, Dv_0)\nonumber\\
&=&-\frac{1}{2} (v_s-v_0,D(v_s-v_0))\nonumber\\
&=&-\frac{1}{2} w^T  E (Dw).\nonumber
\end{eqnarray}

To conclude, the energy difference becomes
\begin{equation}\label{energydiff}
    \begin{aligned}
    \mathcal{E}_{n+1}-\mathcal{E}_n\leq&-\frac{1}{2} w^T E (Dw)+\frac{L}{2}  w^T w + w^T \left(\frac{1}{\tau} \hat{B} E_L (G^{-1}w)
    +\hat{B} A E_L (Dw)\right)\\
    &+\alpha w^T  Q (Dw)-\beta w^T Q (w)\\
    =& \frac{1}{\tau} w^T H_0 (G^{-1}w)-w^T H_1 w + w^T H_2 (Dw),
    \end{aligned}
\end{equation}
where $H_0, H_1$ and $H_2$ are given in (\ref{H0} - \ref{H2}). 

Now if $H_0$, $H_1$ and $H_2$ are positive-definite, then the energy difference on the left-hand-side of (\ref{energydiff}) is negative. It is easy to verify that if $Q$ and $H_0$ are positive-definite, then $H_1$ and $H_2$ are also positive-definite provided that (\ref{alpbet}) is satisfied.

\end{proof}

\begin{remark}
We remark Theorem \ref{maintheorem} that if the condition of positive definiteness on the three matrices is satisfied, then the IMEX-RK methods can decrease the energy unconditionally. We point out that the condition is only a sufficient one, which is confirmed by our numerical experiments, see Example \ref{ex6.2}.
\end{remark}

\begin{remark}
In fact the splitting $\frac{\alpha}{2} |D_{1/2} u|^2 $ is unnecessary for certain Runge--Kutta methods. We may notice that $H_2(0)>0$ guarantees the positive-definiteness of $H_2(\alpha)$ when $Q$ is positive-definite. However, in order to derive unconditional energy dissipation, the stabilization term $\frac{\beta}{2} |u|^2$ is necessary. It is also interesting to notice that the stabilization can improve the numerical performance significantly, see \cite{SAV2}. 
\end{remark}

\begin{remark}
For the MBE model without slope selection, the proof is slightly different but shares the same key idea as long as we consider $w=(\nabla(v_1-v_0),...,\nabla(v_s-v_{s-1}))^T$ instead and replace the key inequality (3.18) in the analysis by the discrete version of (2.3). We omit the detailed proof to keep the presentation short.
\end{remark}

\begin{remark}
In (\ref{energydiff}), the matrices $H_0, H_1$ and $ H_2$ can be viewed as coefficients or some combination of $H^{-k},L^2$ and $H^1$ norms of the difference of the solution at different time steps, where $k$ depends on the operator $G$. Therefore, we do not have to require all three terms to be positive-definite, but could use $H^{-k},H^1$ terms to cover the $L^2$ term.
\end{remark}

\subsection{Application to Allen--Cahn and Cahn--Hilliard equations}

In the remaining of this section, we make a direct application of Theorem \ref{maintheorem} to the Allen--Cahn and the Cahn--Hilliard equations.

\begin{theorem}\label{ACtheorem}
(For the Allen--Cahn equation) For the IMEX-RK scheme (\ref{IMEXRKsys}), if $H_0= (\hat{A})^{-1} E_L$ and $Q=\left((\hat{A})^{-1} A-I\right)E_L+I$ are both positive-definite, then it decreases the energy of the Allen--Cahn equation (\ref{AC}) when following inequalities hold:
\begin{eqnarray}
    && \frac{1}{\tau} \lambda_{\min}(H_0)+\beta \lambda_{\min}(Q) \geq \frac{L}{2} \\
    && \alpha \lambda_{\min}(Q)\geq -\lambda_{\min}((\hat{A})^{-1} A E_L -\frac{1}{2} E ).
\end{eqnarray}

\end{theorem}

\begin{theorem}
(For the Cahn--Hilliard equation) For the IMEX-RK scheme (\ref{IMEXRKsys}), if $H_0= (\hat{A})^{-1} E_L$ and $Q=\left((\hat{A})^{-1} A-I\right)E_L+I$ are both positive-definite, then it decreases the energy of the Cahn--Hilliard equation (\ref{CH}) when the following inequality holds:
\begin{equation}
    \frac{4}{\tau} \lambda_{\min}(H_0) \lambda_{\min}(H_2)\geq (\frac{L}{2}-\beta \lambda_{\min}(Q) )^2
\end{equation}
where $H_2=\alpha Q+(\hat{A})^{-1} A E_L-\frac{1}{2} E $.
\end{theorem}
\begin{proof}
Using $\|(D_{1/2})^{-1} u\|_2 \|D_{1/2} u\|_2\geq \|u\|_2^2$ to simplify (3.22) leads to the result. More precisely,
\begin{equation}
\begin{aligned}
    &-\frac{1}{\tau} w^T H_0 (G^{-1}w) - w^T H_2 (Dw) \\
    \geq &\frac{1}{\tau} \lambda_{\min}(H_0) \| (D_{1/2})^{-1} w\|_2^2 + \lambda_{\min}(H_2)\|D_{1/2} w\|_2^2 \\
    \geq & 2\sqrt{\frac{1}{\tau}\lambda_{\min}(H_0) \lambda_{\min}(H_2)} \|w\|_2^2 \\
    \geq & |\frac{L}{2}-\beta \lambda_{\min}(Q)| \|w\|_2^2 \\
    \geq & -w^T H_1 w,
\end{aligned}
\end{equation}
and thus the energy is decreasing.

\end{proof}

\section{Error Analysis}\label{sec4}

In this section, an error analysis based on the truncation error will be provided. We assume that the exact solutions are sufficiently smooth. For simplicity, we only study the Allen--Cahn and the Cahn--Hilliard models, which are representative $L^2$ and $H^{-1}$ gradient flows. Hence the operator $D$ is fixed as $\Delta$ with $G=-1$ for the Allen--Cahn equations and $G=\Delta$ for the Cahn--Hilliard equations. We point out that the extension to other models would be straightforward as long as the operators $D$ and $G$ induce inner products. We present a framework for $p$th-order IMEX-RK schemes which has $O(\tau^{p+1})$ truncation error, and for more details concerning the order conditions, see Remark \ref{rmk42} and \cite{Ascher,HighlyStabIMEXRK} for example.

\begin{theorem}\label{errthm}
Consider the Allen--Cahn type and the Cahn--Hilliard type phase-field models. Assume the IMEX-RK scheme (\ref{IMEXRKsys}) has $O(\tau^{p+1})$ truncation error, the exact solution $u(t)$ is sufficiently smooth, and the conditions \eqref{H0}-\eqref{alpbet} in Theorem \ref{maintheorem} are satisfied, the following estimate for the IMEX-RK scheme holds:
\begin{equation}\label{errestimate}
        \Vert  u(t_n)-u_n \Vert_{H^1} \leq C e^{C T} \tau^{p},
\end{equation}
where $C$ is a positive constant dependent on the Lipchitz constant, the stabilizer constant, the smoothness of exact solutions but independent of $\tau$.
\end{theorem}
\begin{proof}
At each stage time step $t_{ni}=t_n+c_i \tau, 1\leq i \leq s$, we define intermediate functions $\Bar{v}_i$ by
\begin{eqnarray}
\begin{aligned}
    \Bar{v}_0&&=&\quad u(t_n),\\
    \begin{pmatrix}
    \Bar{v}_1 \\ \Bar{v}_2 \\ \vdots \\ \Bar{v}_s
    \end{pmatrix}
    &&=&\quad
    \begin{pmatrix}\label{refersol}
    \Bar{v}_0 \\ \Bar{v}_0 \\ \vdots \\ \Bar{v}_0
    \end{pmatrix}
    +\tau\left(A \begin{pmatrix}
    \mathcal{L}\Bar{v}_1 \\ \mathcal{L}\Bar{v}_2 \\ \vdots \\ \mathcal{L}\Bar{v}_s
    \end{pmatrix}+\hat{A} \begin{pmatrix}
    N(\Bar{v}_0) \\ N(\Bar{v}_1) \\ \vdots \\ N(\Bar{v}_{s-1})
    \end{pmatrix}    \right),
   \end{aligned}
\end{eqnarray}
and we compute the local truncation error $r_{n+1}$ by plugging the exact solution and the intermediate functions $\Bar{v}_i$  into the scheme as
\begin{equation}
    u(t_{n+1})= u(t_n)
    +\tau \left(b^T \begin{pmatrix}\label{localte}
    \mathcal{L}\Bar{v}_1 \\ \mathcal{L}\Bar{v}_2 \\ \vdots \\ \mathcal{L}\Bar{v}_s
    \end{pmatrix} 
    +\hat{b}^T \begin{pmatrix}
    N(\Bar{v}_0) \\ N(\Bar{v}_1) \\ \vdots \\ N(\Bar{v}_{s-1})
    \end{pmatrix}\right)+r_{n+1}.
\end{equation}
The truncation error has been analyzed in many related papers, see e.g., \cite{Ascher}. More precisely, considering the Taylor expansion of (\ref{refersol})-\eqref{localte} at $t=t_n$ and using asymptotic analysis and the order conditions \eqref{odc12},
we can derive the bound of the local truncation error/ residual as $\|r_{n+1}\|_{H^1}\le C \tau^{p+1}.$

In order to obtain error estimates, we define $e_j=\Bar{v}_j-v_j, 0\leq j\leq s-1$ and $e_s=u(t_{n+1})-u_{n+1}=\Bar{v}_s-v_s-r_{n+1}$. The difference of (\ref{refersol}) and (\ref{IMEXRKsysalter}) shows
\begin{equation}
    \begin{pmatrix}
    e_1 \\ e_2 \\ \vdots \\ e_s
    \end{pmatrix}
    =
    \begin{pmatrix}
    e_0 \\ e_0 \\ \vdots \\ e_0
    \end{pmatrix}
    +\tau\left(A \begin{pmatrix}
    \mathcal{L}e_1 \\ \mathcal{L}e_2 \\ \vdots \\ \mathcal{L}e_s
    \end{pmatrix}+\hat{A} \begin{pmatrix}
    N(\Bar{v}_0)-N(v_0) \\ N(\Bar{v}_1)-N(v_1) \\ \vdots \\ N(\Bar{v}_{s-1})-N(v_{s-1})
    \end{pmatrix}    \right)
    +\begin{pmatrix}
    0 \\ 0 \\ \vdots \\ r_{n+1}
    \end{pmatrix}.
\end{equation}
Using $ \mathcal{L}v=G D_s v,N(v)=-G f_s(v)$, where $D_s=-(1+\alpha)D+\beta I$ and $f_s=-f-\alpha D+\beta I$, recalling that for both Allen--Cahn and Cahn--Hilliard models $D=\Delta$, and following all steps in the proof of Theorem \ref{maintheorem}, we derive
\begin{equation}\label{errineq1}
\begin{aligned}
 &\frac{1}{2}\Vert \nabla e_s\Vert_2^2-\frac{1}{2}\Vert \nabla e_0\Vert_2^2 \\
 =& \frac{1}{\tau} q^T H_0 (G^{-1} q)-\beta q^T Q q + q^T H_2 (Dq)
    - q^T \begin{pmatrix}
    f(\Bar{v}_0)-f(v_0) \\ f(\Bar{v}_1)-f(v_1) \\ \vdots \\ f(\Bar{v}_{s-1})-f(v_{s-1})
    \end{pmatrix}-\frac{1}{\tau} q^T (\hat{A})^{-1} \begin{pmatrix}
    0 \\ 0 \\ \vdots \\ r_{n+1}
    \end{pmatrix},
\end{aligned}
\end{equation}
where $q=(e_1-e_0,e_2-e_1,\cdots,e_s-e_{s-1})^T$. Under the assumptions in Theorem \ref{maintheorem}, $H_0,Q$ and $H_2$ are all positive-definite, so the first three terms on the right-hand side are all negative. For the other two terms, we have
\begin{eqnarray}\label{errstep1}
        -\frac{1}{\tau} q^T (\hat{A})^{-1} \begin{pmatrix}
    0 \\ 0 \\ \vdots \\ r_{n+1}
    \end{pmatrix}&\leq& \frac{C_1}{\tau} \Vert r_{n+1} \Vert_2^2+\frac{C_2}{\tau} q^T q,\\
    -q^T \begin{pmatrix}
    f(\Bar{v}_0)-f(v_0) \\ f(\Bar{v}_1)-f(v_1) \\ \vdots \\ f(\Bar{v}_{s-1})-f(v_{s-1})
    \end{pmatrix} 
    &\leq& \frac{C_3}{\tau} q^T q+C_4 \tau \sum_{i=0}^{s-1} \|e_i\|^2 \label{errstep2} \\
    &\leq& \frac{C_3}{\tau} q^T q+C_4 \tau \sum_{i=0}^{s-1} \|\sum_{j=1}^i (e_j-e_{j-1})+e_0\|^2\nn  \\
    &\leq& \frac{C_3}{\tau} q^T q+C_4 \tau \sum_{i=0}^{s-1} (\sum_{j=1}^i \|e_j-e_{j-1}\|^2+\|e_0\|^2)\nn  \\
    &\leq& \frac{C_3}{\tau} q^T q+C_5 \tau q^T q + C_6 \tau \Vert e_0 \Vert_2^2\nn  \\
    &\leq& \frac{C_3}{\tau} q^T q+C_5 \tau q^T q + C_6 C(\Omega) \tau \Vert \nabla e_0 \Vert_2^2,\nn
\end{eqnarray}
where $C_1, C_3$ can be arbitrary constants and will be determined later, $C_2=C(\hat{A})/C_1, C_4=C(L)/C_3$ and $C_5=C_6=s C_4$.

For the Allen--Cahn equation, adding up all inequalities above (\ref{errineq1}-\ref{errstep2}) and setting $C_1+C_3\leq \lambda_{min}(H_0)$ and $C_5\tau\leq \beta \lambda_{min}(Q)$ lead to the following result
\begin{equation}\label{errest}
    \Vert \nabla e_s \Vert_2^2 \leq (1+C\tau)\Vert \nabla e_0 \Vert_2^2  +\frac{C'}{\tau} \Vert r_{n+1} \Vert_2^2,
\end{equation}
where $e_s=u(t_{n+1})-u_{n+1}$ and $e_0=u(t_n)-u_n$. For the Cahn--Hilliard equation, $q_i$'s are all zero mean functions and for all $u, \ v$ we have
\begin{equation}
    (u,v)\leq \|u\|_{H^1} \|v\|_{H^{-1}} .
\end{equation}
Similarly we have the following results
\begin{eqnarray}\label{errstep3}
        -\frac{1}{\tau} q^T (\hat{A})^{-1} \begin{pmatrix}
    0 \\ 0 \\ \vdots \\ r_{n+1}
    \end{pmatrix}&\leq& \frac{C_1}{\tau} \Vert \nabla r_{n+1} \Vert_2^2+\frac{C_2}{\tau} \Vert (\nabla)^{-1} q \Vert_2^2,\\
    -q^T \begin{pmatrix}
    f(\Bar{v}_0)-f(v_0) \\ f(\Bar{v}_1)-f(v_1) \\ \vdots \\ f(\Bar{v}_{s-1})-f(v_{s-1})
    \end{pmatrix} 
    &\leq& \frac{C_3}{\tau} \Vert (\nabla)^{-1} q \Vert_2^2+C_5 \tau \Vert \nabla q \Vert_2^2 + C_6 \tau \Vert\nabla e_0 \Vert_2^2.\label{errstep4}  
\end{eqnarray}
Thus we could derive an inequality for the Cahn--Hilliard equation
\begin{equation}\label{errest2}
    \Vert \nabla e_s \Vert_2^2 \leq (1+C\tau)\Vert \nabla e_0 \Vert_2^2  +\frac{C'}{\tau} \Vert \nabla r_{n+1} \Vert_2^2.
\end{equation}

Therefore, by Gronwall's inequality, we derive the error estimate (\ref{errestimate}).
\end{proof}
\begin{remark}
    In the proof, $\Vert (\nabla^{-1}) u \Vert_2=-(u,D^{-1} u)$ is well-defined when considering the Cahn--Hilliard type equations, where $u$ has zero mean.
\end{remark}
\begin{remark}\label{rmk42}
    The IMEX-RK scheme (\ref{IMEXRKsys}) is of order $p$ if the
following order conditions are satisfied for $1\leq j \leq p\leq 4$
 \begin{eqnarray}
 \begin{aligned}
      & b_\sigma^T c_\sigma^{j-1}=\hat{b}_\sigma^T c_\sigma^{j-1}=\frac{1}{j}; \quad
      b_\sigma^T \Bar{A}_\sigma^{j-1} e_\sigma =\hat{b}_\sigma^T \Bar{A}_\sigma^{j-1} e_\sigma =\frac{1}{j!}; \\
      & (b_\sigma\cdot c_\sigma)^T \Bar{A}_\sigma c_\sigma=(\hat{b}_\sigma\cdot c_\sigma)^T \Bar{A}_\sigma c_\sigma=\frac{1}{8}, \quad
       b_\sigma^T \Bar{A}_\sigma  c_\sigma^2 =\hat{b}_\sigma^T \Bar{A}_\sigma c_\sigma^2 =\frac{1}{12},\ \text{if $p=4$},
      \end{aligned}\label{odc12}
\end{eqnarray}
where $\sigma=s+1,c_\sigma^{j}=(0,c_1^j,c_2^j,...,c_s^j)^T$, each $\Bar{A}_\sigma$ can be either $A_\sigma$ or $\hat{A}_\sigma$, and $b_\sigma,\hat{b}_\sigma,A_\sigma,\hat{A}_\sigma$ represent $\sigma\times1$ vectors and $\sigma\times\sigma$ matrices in (\ref{IMEXRK}) respectively.
\end{remark}

\section{Implicit-Explicit Runge-Kutta schemes}\label{sec5}

In this section we present some Runge--Kutta schemes to illustrate how the conditions of the theorem are satisfied and offer a new third-order four-stage IMEX-RK scheme which unconditionally decrease the energy of the gradient flows under our consideration.

\subsection{Example 1: first-order IMEX}
This simple one-step case corresponds to the following IMEX scheme whose tableau reads:
\begin{equation}
\begin{array}{c|cc}
    0   & 0 & 0 \\
    1   & 0 & 1 \\
    \hline
        & 0 & 1
\end{array}, \quad 
\begin{array}{c|cc}
    0   & 0 & 0 \\
    1   & 1 & 0 \\
    \hline
        & 1 & 0
\end{array}.
\end{equation}
Here $A=\hat{A}=1$, and thus $H_0=Q=1$ and $H_2(0)=1/2$. Therefore, for the Allen-Cahn equation, we only need
\begin{equation}
    \frac{1}{\tau}+\beta \geq 1
\end{equation}
to guarantee the energy dissipation law; while for the Cahn--Hilliard equation, we only need
\begin{equation}
    \frac{2\epsilon^2}{\tau}+\beta \geq \frac{L}{2}.
\end{equation}

\subsection{Example 2: a second-order IMEX}

Consider this second-order IMEX-RK scheme with coefficients (see, e.g.,\cite{Ascher})
\begin{equation}
\begin{array}{c|ccc}
    0      & 0 & 0 & 0 \\
    \gamma & 0 & \gamma & 0 \\
    1      & 0 & 1-\gamma &\gamma\\
    \hline
           & 0 & 1-\gamma &\gamma
\end{array}, \quad 
\begin{array}{c|ccc}
    0      & 0 & 0 & 0 \\
    \gamma & \gamma & 0 & 0 \\
    1      & \delta  & 1-\delta &0\\
    \hline
           & \delta  & 1-\delta &0
\end{array},
\end{equation}
where $\gamma=1-\frac{\sqrt{2}}{2}, \delta=1-\frac{1}{2\gamma}$. Here
\begin{equation}
    A=\begin{pmatrix}
    \gamma & 0 \\ 1-\gamma &\gamma
    \end{pmatrix},\quad \hat{A}=\begin{pmatrix}
    \gamma & 0 \\ \delta  & 1-\delta
    \end{pmatrix},
\end{equation}
and thus 
\begin{equation}
    H_0=\begin{pmatrix}
    2+\sqrt{2} & 0 \\ 2 & 2-\sqrt{2}
    \end{pmatrix},\quad
    Q=\begin{pmatrix}
    1 & 0\\0 & 3-2\sqrt{2}
    \end{pmatrix},\quad
    H_2=\begin{pmatrix}
    1/2 & -1/2 \\ 1/2 & 5/2-2\sqrt{2}
    \end{pmatrix}.
\end{equation}
The corresponding smallest eigenvalues of their symmetrizers are
\begin{equation}
    \lambda_{\min}(H_0)=2-\sqrt{3},\quad \lambda_{\min}(Q)=3-2\sqrt{2},\quad \lambda_{\min}(H_2(0))=5/2-2\sqrt{2}.
\end{equation}
Therefore, we need the $\alpha |D_{1/2} u|^2$ term in the splitting. For the Allen-Cahn equation, we need
\begin{equation}
    \begin{aligned}
    & \frac{2-\sqrt{3}}{\tau} + (3-2\sqrt{2}) \beta \geq 1,\\
    & \alpha \geq \alpha_0=\frac{2\sqrt{2}-5/2}{3-2\sqrt{2}}\approx 1.9142.
    \end{aligned}
\end{equation}
For the Cahn--Hilliard equation, we need  
\begin{equation}
    4(2-\sqrt{3})(\alpha +5/2-2\sqrt{2})\frac{\epsilon^2}{\tau}+(3-2\sqrt{2})\beta \geq \frac{L}{2}.
\end{equation}

However, since $\lambda_{\min}(Q)$ for this scheme is too small, if we want to get unconditional energy dissipation we have to set $\beta>\frac{L}{2\lambda_{\min}(Q)}=3+2\sqrt{2}$, which might be too large and may cause larger truncation errors.

\subsection{Third-order schemes}
There is no pair of a three-stage, $L$-stable DIRK (diagonally implicit Runge--Kutta) and a four-stage ERK (explicit Runge--Kutta) with a combined third-order accuracy (see \cite{Ascher}), so we have to consider four-stage schemes. However, as far as we search, no existing four-stage third-order ARS Runge--Kutta scheme satisfies the conditions in Theorem \ref{maintheorem}. This motivates us to construct a new one. Below we first list two common four-stage third-order Runge--Kutta schemes as examples to illustrate why they fail, and then present our new four-stage third-order scheme.

\vspace{.5em}
\textit{Example 4:}
\begin{equation}
    A=\begin{pmatrix}
    1/2  & 0   & 0 & 0 \\ 
    1/6  & 1/2 & 0 & 0 \\
    -1/2 & 1/2 & 1/2 & 0\\
    3/2 & -3/2 & 1/2 & 1/2
    \end{pmatrix},\quad \hat{A}=\begin{pmatrix}
    1/2  & 0 & 0 & 0 \\ 
    11/18 & 1/18 & 0 & 0 \\
    5/6  & -5/6 & 1/2 & 0\\
    1/4 & 7/4 & 3/4 & -7/4
    \end{pmatrix}.
\end{equation}
The corresponding smallest eigenvalues of symmetrizers of the determinants are approximately
\begin{equation}
    \lambda_{\min}(H_0)=-15.242727,\quad \lambda_{\min}(Q)=-7.706226,\quad \lambda_{\min}(H_2(0))=-8.206226.
\end{equation}
Since both $H_0$ and $Q$ are negative-definite, this scheme also does not satisfy the conditions of our theorem.

\vspace{.5em}
\textit{Example 5: Energy decreasing four-stage third-order IMEX-RK scheme.} Here we present a four-stage third-order ARS IMEX-RK scheme which has rational coefficients in $c$ (time step of stages) which are not too large.

\begin{equation}
\begin{array}{c|ccccc}
    0      & 0 & 0 & 0 & 0 & 0 \\
    3/5    & 0 & 0.6 & 0 & 0 & 0\\
    3/2    & 0 & 0.46875 & 1.03125 & 0 & 0\\
    19/20  & 0 & 0.4 & -0.5578125 & 1.1078125 & 0\\
    1      & 0 & a_{41} & a_{42} & a_{43} & 25.75\\
    \hline
           & 0 & a_{41} & a_{42} & a_{43} & 25.75
\end{array}
\end{equation}
\begin{equation}
\begin{array}{c|ccccc}
    0      & 0 & 0 & 0 & 0 & 0 \\
    3/5    &  0.6 & 0 & 0 & 0 & 0\\
    3/2    &  0.796875 & 0.703125 & 0 & 0 & 0\\
    19/20  &  0.4 & \hat{a}_{32} & \hat{a}_{33} & 0 & 0\\
    1      & \hat{a}_{41} & \hat{a}_{42} & \hat{a}_{43} & \hat{a}_{44} & 0\\
    \hline
           & \hat{a}_{41} & \hat{a}_{42} & \hat{a}_{43} & \hat{a}_{44} & 0
\end{array}
\end{equation}
where
\begin{equation}
    \begin{aligned}
    a_{41}=3.736772486772523;\quad\quad&a_{42}=-0.781144781144795;\\
    a_{43}=-27.705627705628103;\;\;\; &\hat{a}_{32}=0.420225694444444;\\
    \hat{a}_{33}=0.129774305555556;\quad\quad&\hat{a}_{41}=0.301169590643275;\\
    \hat{a}_{42}=0.330687830687831;\quad\quad&\hat{a}_{43}=-0.087542087542087;\\
    \hat{a}_{44}=0.455684666210982;\quad\quad&
    \end{aligned}
\end{equation}
The corresponding smallest eigenvalues of symmetrizers of the determinants are approximately
\begin{equation}
    \lambda_{\min}(H_0)\approx0.087230,\quad \lambda_{\min}(Q)=1,\quad \lambda_{\min}(H_2(0))=0.5,
\end{equation}
which are all positive. Therefore, we only need to set 
\begin{equation}
    \beta\geq\frac{1}{\lambda_{\min}(Q)}=1,
\end{equation}
so that this scheme unconditionally decreases the energy of phase fields with Lipschitz nonlinear terms.

Here we briefly illustrate a strategy to construct four-stage third-order IMEX-RK schemes whose coefficients satisfy the conditions in Theorem 3.1. Since we only consider type ARS IMEX-RK schemes, we automatically have
\begin{equation}
c_\sigma=\hat{c}_\sigma,\;  b_\sigma=A_\sigma(\sigma,:), \;
 \hat{b}_\sigma=\hat{A}_\sigma(\sigma,:).
\end{equation}
Therefore, we have 24 unknowns in total, including 4 in $c_\sigma$, 10 in $A_\sigma$ and 10 in
$\hat{A}_\sigma$. The order conditions for IMEX-RK3 are
\begin{eqnarray}
      & A_\sigma e_\sigma = \hat{A}_\sigma e_\sigma =c_\sigma,
      \\  \label{vdmcond}
      & b_\sigma^T e_\sigma=\hat{b}_\sigma^T e_\sigma=1,\ b_\sigma^T c_\sigma=\hat{b}_\sigma^T c_\sigma=\frac{1}{2},\  b_\sigma^T c_\sigma^2=\hat{b}_\sigma^T c_\sigma^2=\frac{1}{3}, \\ \label{Asol}
      & b_\sigma^T A_\sigma c_\sigma =\hat{b}_\sigma^T A_\sigma c_\sigma =b_\sigma^T \hat{A}_\sigma c_\sigma =\hat{b}_\sigma^T \hat{A}_\sigma c_\sigma =\frac{1}{6}.
\end{eqnarray}
Notice that the first two equations in (5.19) are the same as the last line of (5.18), thus we have 16 different equations in total. Although the whole system is nonlinear, the relation between these variables can be simplified. In particular, by setting some free parameters, we can find a linear way to search the coefficients. 

First we set $c_\sigma=(0, c_1, c_2, c_3, c_4)$, where  $c_4=1$ can be derived by the stiffly accurate condition and $b_\sigma^T e_\sigma =1$ and $c_1,c_2, c_3$ are different free parameters (also not equal to 1).

Then we set $b_\sigma^T c_\sigma^3=\xi, \hat{b}_\sigma^T c_\sigma^3=\hat{\xi}$ where
$\xi$ and $\hat{\xi}$ are also free parameters. Now we could solve $b_\sigma$ and $\hat{b}_\sigma$ (depending on $c_\sigma$) by 
\begin{equation}
b_\sigma^T \left(e_\sigma, c_\sigma, c_\sigma^2, c_\sigma^3\right) =(1,1/2,1/3,\xi),\quad  \hat{b}_\sigma^T \left(e_\sigma, c_\sigma, c_\sigma^2, c_\sigma^3\right) =(1,1/2,1/3,\hat{\xi}),
\end{equation}
where the Vandermonde matrix is nondegenerate since $c_1, c_2, c_3, c_4$
are different. 

Finally, with the knowledge of $c_\sigma$, $b_\sigma$ and $\hat{b}_\sigma$, the only nonlinear equation (\ref{Asol}) becomes a simple linear problem for $A_\sigma$ and $\hat{A}_\sigma$ which is easy to solve.

Therefore, we can set the range for all these free parameters and find all possible combinations of them. By doing this we can make the Runge-Kutta tables satisfy the conditions of Theorem 3.1.

However, such trick does not work for fourth and higher order situations, since the corresponding order condition of (\ref{Asol}) will have forms
\begin{equation}
    b_\sigma^T A_\sigma^{p-2} c_\sigma =\hat{b}_\sigma^T A_\sigma^{p-2} c_\sigma =b_\sigma^T \hat{A}_\sigma^{p-2} c_\sigma =\hat{b}_\sigma^T \hat{A}_\sigma^{p-2} c_\sigma =\frac{1}{p!}.
\end{equation}
Note that even if we already know all variables in $b$ and $c$, a nonlinear problem still needs to be solved. Therefore, it is still difficult to find a fourth-order scheme.

\section{Numerical Tests}\label{sec6}

In this section numerical experiments will be carried out to confirm our theoretical results. In all examples, we assume the periodic boundary conditions, fix attention on the $[0,2\pi]$ domain, and use a Fourier spectral method for space variables. Besides, we also consider a truncated double-well potential $\Tilde{F}(u)$ so that it naturally satisfies the Lipschitz condition. More precisely, for a sufficiently large $M$ ($M=2$ is enough for cases in this paper), we replace $F(u)=\frac{1}{4} (u^2-1)^2$ by 
\begin{eqnarray}\nn
    \Tilde{F}(u) &= \left\{\begin{aligned}
        &\frac{3M^2-1}{2} u^2 - 2 \text{sgn}(u)  M^3 u +\frac{1}{4} (3M^4+1), \quad &|u|>M  \\
        &\frac{1}{4} (u^2-1)^2, \quad &|u|\leq M  
    \end{aligned}
    \right.,
\end{eqnarray}
and $f(u)=u^3-u$ by
\begin{eqnarray}\nn
    \Tilde{f}(u)=\Tilde{F}'(u) &= \left\{\begin{aligned}
        &(3M^2-1) u - 2 \text{sgn}(u)  M^3, \quad &|u|>M  \\
        &u^3-u, \quad &|u|\leq M  
    \end{aligned}
    \right..
\end{eqnarray} In fact, the maximum norm of numerical solutions never exceeds the bound $M$ so this replacement does not affect the properties of numerical solutions.

\subsection{Accuracy Test}
\begin{figure}[h]
    \centering
    \includegraphics[width=7cm,height=5cm]{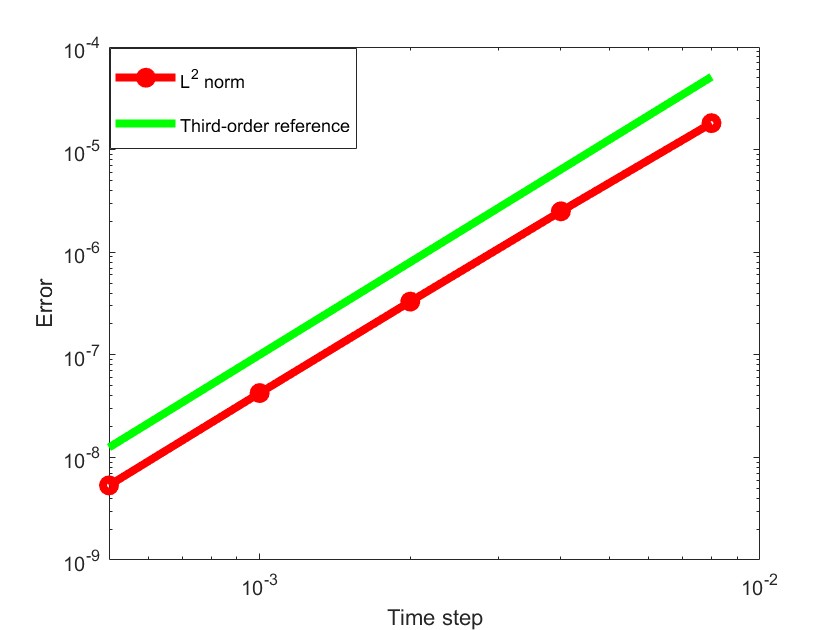}
    \includegraphics[width=7cm,height=5cm]{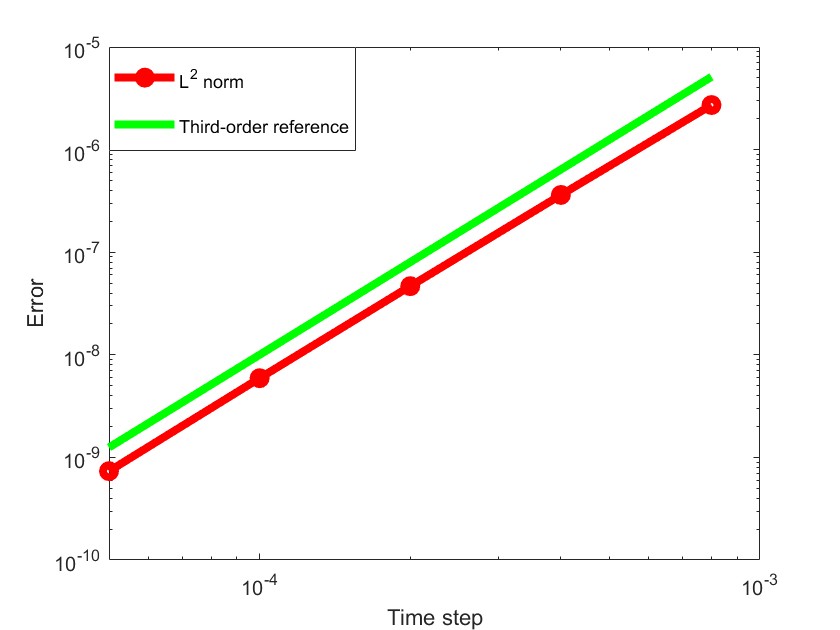}
    \caption{Accuracy Tests for Allen--Cahn (left) and Cahn--Hilliard (right) models}
    \label{FigAccTest}\vspace{-1em}
\end{figure}

 We take the energy-decreasing four-stage third-order IMEX-RK scheme as the example to test the accuracy. We apply this IMEX-RK method for the time discretization to solve the Allen--Cahn equation (\ref{AC}) and the Cahn--Hilliard equation (\ref{CH}). The initial condition is given by $u_0=0.05\sin(x)\sin(y)$. The number of discrete points in space is $N=128$ and the parameter $\epsilon=0.1$ for the Allen--Cahn equation and $\epsilon=1$ for the Cahn--Hilliard equation. The reference solution is given by this IMEX-RK method with sufficiently small time step $\tau=10^{-5}$ (we denote '$dt$' as time step in all figures). The numerical errors, represented by $L^2$ norms at $T=0.032$, and the corresponding rates are shown in Table \ref{tabACerr} and \ref{tabCHerr}. We also plot the relationship between numerical errors at a longer time scale $T=0.128$ and time step sizes in Fig \ref{FigAccTest}. From both numerical experiments we can observe the third order convergence for both equations.
\begin{table}[h]
     \centering
     \begin{tabular}{|c|c|c|c|c|c|}
        \hline
         Allen--Cahn & $ \tau=8e-3$ & $ \tau=6.4e-3$ & $ \tau=4e-3$ & $ \tau=2e-3$   & $ \tau=1e-3$       \\
         \hline
         $L^2$ Error & 4.123e-06 & 2.195e-06 & 5.698e-07 & 7.513e-08  & 9.647e-09\\
         \hline
         Rate & - & 2.855 & 2.923  & 2.961 & 2.990    \\
        \hline
        \end{tabular}
        \caption{Errors and rates for the Allen--Cahn equation}
     \label{tabACerr}
     \centering
     \begin{tabular}{|c|c|c|c|c|c|}
        \hline
         Cahn--Hilliard & $ \tau=8e-4$ & $ \tau=4e-4$ & $ \tau=2e-4$ & $ \tau=1e-4$  & $ \tau=5e-5$ \\
         \hline
         $L^2$ Error & 1.506e-08 & 7.846e-09 &  1.979e-09 & 1.026e-9 & 2.543e-10 \\
         \hline
         Rate & - & 2.883 & 2.910 & 2.928 & 2.960 \\
        \hline
        \end{tabular}
        \caption{Errors and rates for the Cahn--Hilliard equation}
     \label{tabCHerr}
\end{table} 

\subsection{Energy Evolution}\label{ex6.2}
 
 Now we present the energy evolution of the four-stage third-order IMEX-RK solutions with different values of the parameter in order to study how the stabilization terms affect the energy stability. 
 
\begin{figure}[h]
    \centering
    \includegraphics[height=6cm,width=10cm]{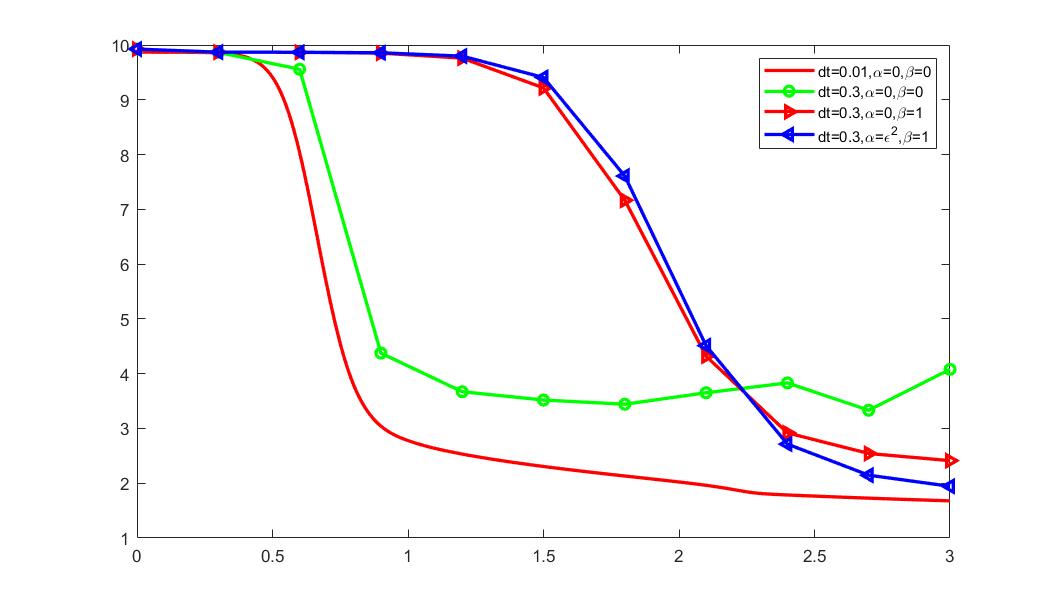}
    \caption{Energy curves of numerical solutions for the Allen--Cahn equation}
    \label{figACenergy}
\end{figure}

 First we take the Allen--Cahn equation (\ref{gf}) with $G=-1/\epsilon, D=\epsilon^2 \Delta$ and $f(u)=u^3-u$ as an example. The scaling of $G$ here is solely for the purpose of better observing numerical experimental phenomena without impacting stability or other properties. It could be viewed as introducing an alternative time scale $t'$ where $dt'=\epsilon dt$, so that the equation is the same as the original model. According to our theorem and detailed analysis in Example $4$, $\alpha=0$ and $\beta=1$ are enough to guarantee the unconditional energy dissipation. We set $N=64$ and $\epsilon=0.1$. Figure \ref{figACenergy} shows that the energy curves of solutions with different parameters starting from the same random initial data.
 
We can observe that when $\alpha=\beta=0$ and $\tau=0.3$, the solution is not accurate enough and the energy may increase, while if we take $\beta=1$, the solution comes back to the energy-decreasing state. Moreover, $\alpha=\epsilon^2$ or $\alpha=0$ does not make much difference in this case.

\begin{figure}[h]
    \centering
    \includegraphics[height=6.5cm,width=10cm]{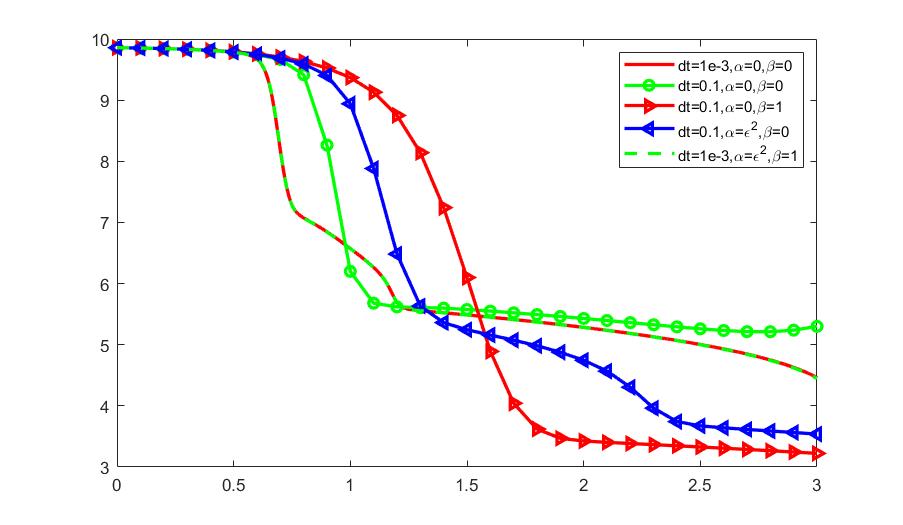}
    \caption{Energy curves of numerical solutions for the Cahn--Hilliard equation}
    \label{figCHenergy}
\end{figure}

Then we take the Cahn--Hilliard equation (\ref{CH}) as an example, starting with the given initial data $u_0=0.05\sin(x)\sin(y)$ since the corresponding solution has its special energy curve. For this case we set $N=128$ and $\epsilon=0.1$. Figure \ref{figCHenergy} shows the development of all these energy curves.

The first and the last curves indicate when the time step is small enough, the solutions of $\alpha=\beta=0$ and $\alpha=\epsilon^2,\beta=1$ cases are quite close, showing the same special evolution of the energy. When $\tau=0.1$ and $\alpha=\beta=0$, the solution seems to be rather inaccurate and the energy may increase. Both $\beta=1$ or $\alpha=\epsilon^2$ cases help the energy to dissipate. In this example, the $\alpha=\epsilon^2$ case happens to produce more accurate pattern compared with the reference solution. This example also indicates that the condition of positive definiteness on the three matrices is not necessary when applied to the Cahn--Hilliard equations.

\subsection{The stabilization effect on the error}

We study the effect of stabilization terms on the error by taking the Allen--Cahn equation (\ref{gf}) with $G=-1/\epsilon, D=\epsilon^2 \Delta$ and $f(u)=u^3-u$ as example. In our computations we set $N=128$, $\epsilon=0.05$ and $T=1$. According to the Theorem 3.2, for the third-order IMEX-RK in Example 5, $\lambda\epsilon/\tau+\beta\geq 1$ guarantees the energy dissipation, where $\lambda\approx0.087230$. Therefore, if the time step is sufficiently small, then we do not need to add a stabilization term. Otherwise, $\beta=1-\lambda\epsilon/\tau$ is required. In order to preserve the energy dissipation, we study how the error changes as the time step varies for $\beta=0$ and $\beta=\max\{1-\lambda\epsilon/\tau,0\}$. 

\begin{figure}[h]
    \centering
    \includegraphics[height=6.5cm,width=10cm]{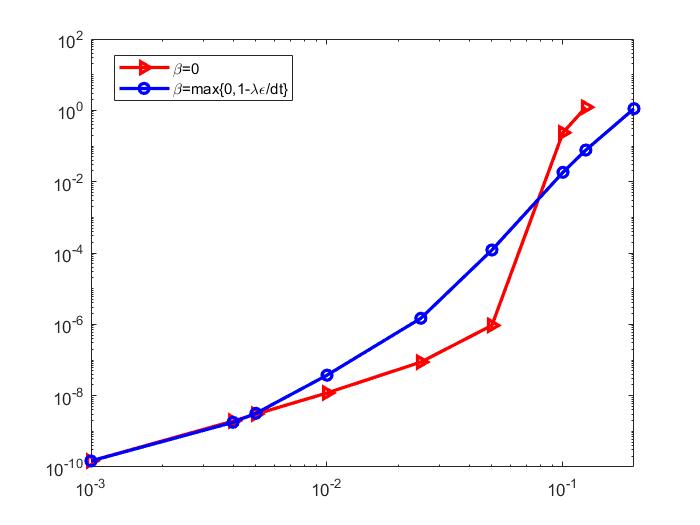}
    \caption{Error-time step curves for the Allen--Cahn equation}
    \label{figACerrdtbeta}
\end{figure}

Figure 3 shows that when the time step is sufficiently small, these two curves overlap. When the time step is big, the $\beta=0$ curve does not show a smooth expansion and the discrete solution with $\tau=0.2$ even blows up. When the time step is small, the $\beta=0$ solution seems to be more accurate. However, if we only need a $10^{-2}$ accuracy, then the solution with $\tau=0.1$ and the stabilization is more accurate than the one without stabilization.

\subsection{The stabilization effect on the dynamics}

In this subsection we show that the stabilization procedure may help to resolve the solution dynamics. We take the Cahn--Hilliard equation (\ref{CH}) as the example, set $N=256$, $\beta=0.1$, $\epsilon=0.02$ and observe the pictures of the solution with $\tau=0.005$ and $\tau=10^{-4}$ and $T=0.1$. 

It is observed from Figure 4 that the large time step solution without stabilization is very inaccurate and has many small points around the boundary and $x=\pi$ and $y=\pi$. From the colorbar we also see that the range of the solution (i.e., the maximal and minimal values) is larger than the results in other three cases. We also present the curves at the section $y=\pi(1+5/N)$ to illustrate what these small points mean in Figure 5.

\begin{figure}[H]
    \centering
    \includegraphics[width=11cm,height=8cm]{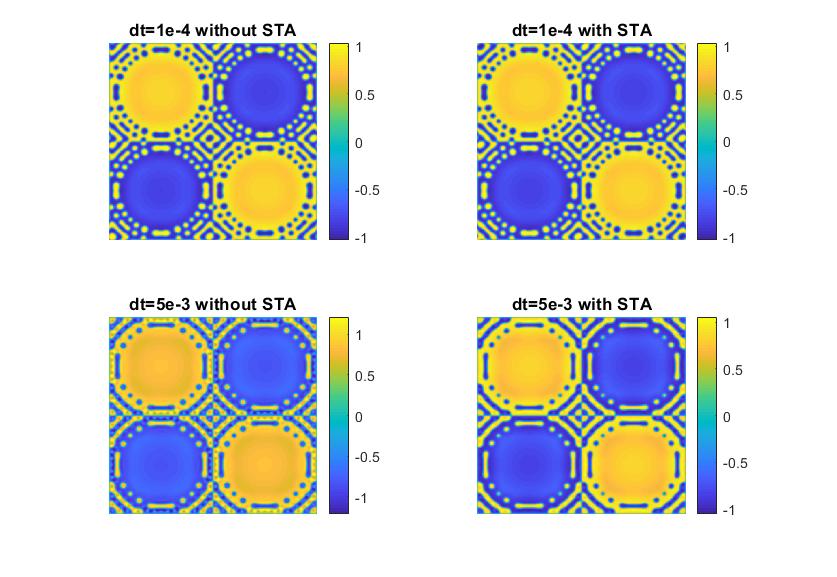}
    \caption{Small and large time step solutions with and without stabilization at $T=0.1$ }
    \label{figACdyna}
\end{figure}

From Figure 5 we observe that the curves obtained without using stabilization have obvious oscillations which are not supposed to happen. 
However, the curves derived by using stabilization seem more reasonable.

\begin{figure}[H]
    \centering
    \includegraphics[width=11cm]{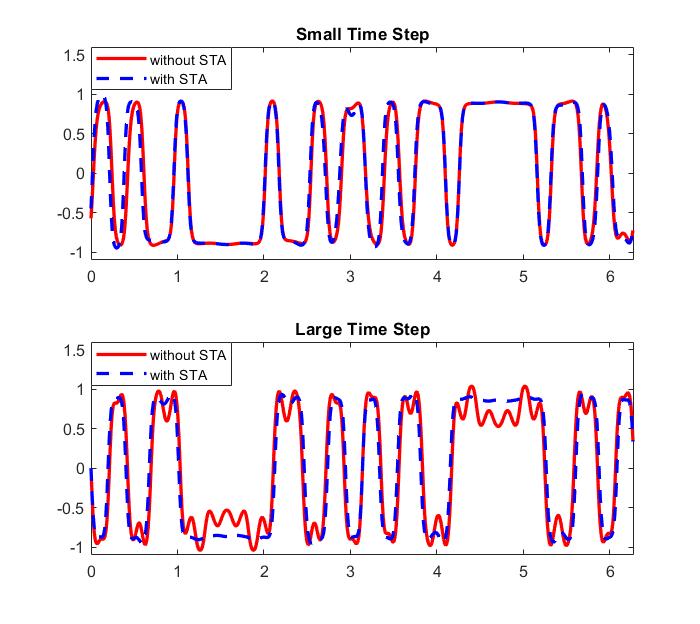}
    \caption{Curves for large and small time step solutions at the section }
    \label{figACdynasec}
\end{figure}

\section{Concluding Remarks}\label{sec7}

In this work, we prove that a class of implicit-explicit Runge--Kutta schemes unconditionally preserve the energy dissipation law for a large family of gradient flows by studying the eigenvalues of matrices which depend on the IMEX-RK tableau. As long as the conditions in our theorems are satisfied, the IMEX-RK 
methods produce energy-diminishing numerical solutions independent of the time step size used. Convergence analysis based on truncation errors is provided. The theoretical predictions are also verified by numerical experiments.



There are two possible directions for future work. One is to construct IMEX-RK of order 4 or higher. It is pointed out at the end of Section 5 that under the current framework, this task seems difficult. Another direction is to extend the IMEX-RK schemes for solving other types of evolution equations satisfying certain stability properties.

\section*{Acknowledgements}

This work is partially supported by the National Science Foundation of China and Hong Kong RGC Joint Research Scheme (NSFC/RGC 11961160718), and the fund of the Guangdong Provincial Key Laboratory of Computational Science and Material Design (No. 2019B030301001). The work of J. Yang is supported by the National Science Foundation of China (NSFC-12271240) and the Shenzhen Natural Science Fund (RCJC20210609103819018). T. Tang is supported in part by the Guangdong Provincial Key Laboratory of Interdisciplinary Research and Application for Data Science under UIC 2022B1212010006.

\bibliography{ref}
\bibliographystyle{siamplain}

\end{document}